\documentclass{svjour3}                     
\smartqed  
\usepackage{graphicx}
\usepackage{amssymb,amsmath,amsfonts,latexsym,euscript}
\def\d{\,{\rm d}}
\def\E{\,{\rm e}}

\def\i{{\rm i}}
\def\abd{\mathop{\ulcorner\!\urcorner}}

\begin{document}

\title{Construction of optimal quadrature formulas exact for exponentional-trigonometric functions by Sobolev's method}

\thanks{A.R. Hayotov is supported by ISEF program at Korea Foundation for Advanced Studies.}

\titlerunning{Optimal formulas for approximate integration}        

\author{A.K.~Boltaev, A.R.~Hayotov$^*$, Kh.M.~Shadimetov}


\institute{A.K. Boltaev \at
V.I.Romanovskiy Institute of Mathematics, Uzbekistan Academy of Sciences,
81 M.Ulugbek str., Tashkent 100170, Uzbekistan\\
\email{aziz\_boltayev@mail.ru}.
\and A.R. Hayotov \at
$^*$Corresponding author:
Department of Mathematical Sciences, KAIST, 291 Daehak-ro, Yuseong-gu, Daejeon 34141, Republic of Korea\\
\email{hayotov@kaist.ac.kr}.\\
V.I.Romanovskiy Institute of Mathematics, Uzbekistan Academy of Sciences,
81 M.Ulugbek str., Tashkent 100170, Uzbekistan,\\
\email{hayotov@mail.ru}.
\and
Kh.M.Shadimetov \at
V.I.Romanovskiy Institute of Mathematics, Uzbekistan Academy of Sciences,
81 M.Ulugbek str., Tashkent 100170, Uzbekistan\\
Tashkent Railway Engineering Institute, 1 Odilxojaev str., Tashkent 100167, Uzbekistan\\
\email{kholmatshadimetov@mail.ru}.
}

\date{Received: date / Accepted: date}

\maketitle
\begin{abstract}
The paper studies Sard's problem on construction of optimal quadrature formulas in the space $W_2^{(m,0)}$ by Sobolev's method.
This problem consists of two parts: first calculating the norm
of the error functional and then finding the minimum of this norm by coefficients of quadrature formulas.
Here the norm of the error functional is calculated with the help of the extremal function. Then using the method of Lagrange multipliers the system of linear equations for coefficients of the optimal quadrature
formulas in the space $W_2^{(m,0)}$ is obtained, moreover the existence and uniqueness of the
solution of this system are discussed. Next, the discrete analogue $D_m(h\beta)$ of the differential operator $\frac{\d^{2m}}{\d x^{2m}}-1$ is constructed.
Further, Sobolev's method of construction of optimal quadrature formulas in the space $W_2^{(m,0)}$, which based on the discrete analogue $D_m(h\beta)$, is described. Finally, for $m=1$ and $m=3$ the optimal quadrature formulas which are exact to exponential-trigonometric functions are obtained.

\textbf{MSC:} 41A05, 41A15, 65D30, 65D32.

\textbf{Keywords:} the extremal function, the error functional, discrete analogue, Sobolev's method, optimal coefficients.
\end{abstract}

\maketitle

\section{Introduction and statement of the problem}

It is known that quadrature formulas are needed for numerical calculation of definite integrals. Additionally, quadrature formulas
provide a basic and important tool for the numerical solution of differential and integral equations.
Development of new algorithms for constructing optimal, in some sense, quadrature formulas as well as an estimation
of their errors in different classes of functions, based on algebraic and variation approaches
is one of the important problems of computational mathematics. The optimization problem of numerical integration formulas
in variation approaches is a problem of finding the minimum of the norm of an error functional on the given space of functions.
There are S.M.Nikol'skii's problem \cite{SMNik50,SMNikBook88} consisting on minimization of the norm of the error functional
by coefficients and by nodes and A.Sard's problem \cite{Sard,MeySar} consisting on minimization of the norm of the
error functional with fixed nodes by coefficients. The solutions of Nikol'skii's and Sard's problems
are called \emph{the optimal quadrature formula} in the sense of Nikol'skii and in the sense of Sard, respectively.

In the present work we study Sard's problem of construction of optimal quadrature formulas in a Hilbert space.

We denote by $W_2^{(m,0)}$ the class  of functions
$\varphi$ defined on the interval $[0,1]$ which posses an absolutely continuous $(m-1)$th  derivative on $[0,1]$
and whose $m$th derivative is in $L_2(0,1)$. The class $W_2^{(m,0)}$, under the pseudo-inner product
\begin{equation}\label{inner-product}
\langle\varphi,\psi\rangle_m=\int\limits_0^1(\varphi^{(m)}(x)+\varphi(x))(\psi^{(m)}(x)+\psi(x))\d x,
\end{equation}
is a Hilbert space if we identify functions that differ by a solution of the equation
$f^{(m)}(x)+f(x)=0$ (see, for example, \cite[p. 213]{Ahlb67}).
Therefore $W_2^{(m,0)}$ is the Hilbert space equipped with the norm
\begin{equation}\label{norm}
\|\varphi\|_{W_2^{(m,0)}}=\left(\int\limits_0^1(\varphi^{(m)}(x)+\varphi(x))^2\d x\right)^{1/2}
\end{equation}
corresponding to the inner product (\ref{inner-product}).

It should be noted that for a linear differential operator of order $m$,
$L\equiv a_m\frac{\d^m}{\d x^m}+a_{m-1}\frac{\d^{m-1}}{\d x^{m-1}}+...+a_1\frac{\d}{\d x}+a_0,$ $a_m\neq 0$,
J.H.~Ahlberg, E.N.~Nilson, and J.L.~Walsh in the book \cite[Chapter 6]{Ahlb67} studied the Hilbert spaces in the context of generalized
splines. Namely, with the pseudo-inner product
\begin{equation*}
\langle\varphi,\psi\rangle_L=\int\limits_a^bL\varphi(x)\cdot L\psi(x)\d x.
\end{equation*}

For a function $\varphi$ from the space $W_2^{(m,0)}$ we consider a quadrature formula of the form
\begin{equation}\label{(1)}
\int\limits_0^1 {\varphi (x)dx} \cong \sum\limits_{\beta  =
0}^N {C_\beta  \varphi (x_{\beta})},
\end{equation}
where $ C_\beta$ are coefficients and $x_\beta$ are nodes situated in the interval $[0,1]$, $\varphi(x_\beta)$ are given values, $N$ is a natural number.
We suppose that $x_{\beta}$ are fixed.

The following difference between integral and quadrature sum
\begin{equation}\label{(4)}
\left(\ell,\varphi\right)=\int\limits_0^1\varphi (x)\d x -\sum\limits_{\beta=0}^N
C_\beta \varphi(x_\beta)
\end{equation}
is called \emph{the error} of the quadrature formula (\ref{(1)}) and
$
(\ell,\varphi)=\int\limits_{-\infty}^{\infty}\ell(x)\varphi(x)\d x
$
is the value of the error functional $\ell$ at the given function $\varphi$. Here the error functional $\ell$
has the form
\begin{equation}\label{(2)}
\ell (x) = \varepsilon _{[0,1]} (x) - \sum\limits_{\beta  =
0}^N {C_\beta  \delta (x - x_\beta  )},
\end{equation}
where $\varepsilon_{[0,1]}(x)$ is the characteristic function of the interval $[0,1]$, $\delta$ is Dirac's delta-function.

According to the Cauchy-Schwarz inequality the absolute value of the error (\ref{(4)}) is estimated by the norm
\begin{equation}\label{normofell}
\left\|{\ell}\right\|_{W_2^{(m,0)*}}=\mathop {\sup }\limits_{\left\| \varphi \right\|_{W_2^{(m,0)}}= 1}
{\left| {\left( {\ell,\varphi}\right)}\right|}
\end{equation}
of the error functional $\ell$ as follows
$$
|(\ell,\varphi)|\leq \|\varphi\|_{W_2^{(m,0)}}\left\|{\ell}\right\|_{W_2^{(m,0)*}}.
$$
Here $W_2^{(m,0)*}$ is the conjugate space to the space $W_2^{(m,0)}$.

Sard's problem on construction of optimal quadrature formulas in the space $W_2^{(m,0)}$
is to find such coefficients $C_{\beta}$ which satisfy the equality
\begin{equation}\label{(5)}
\left\|\mathring{\ell}\right\|_{W_2^{(m,0)*}} = \mathop{\inf}\limits_{C_\beta}
\left\|{\ell}\right\|_{W_2^{(m,0)*}},
\end{equation}
i.e., to find the minimum of the norm (\ref{normofell}) of the error functional $\ell$
by coefficients $C_{\beta}$ for fixed $x_{\beta}$.

This problem consists of two parts: first calculating the norm
(\ref{normofell}) of the error functional $\ell$ in the space $W_2^{(m,0)*}$ and then
finding the minimum of the norm (\ref{normofell}) by coefficients $C_{\beta}$ for fixed nodes $x_{\beta}$.

There are several methods for constructing of optimal quadrature formulas in the
sense of Sard such as the spline method, the $\phi$-function method (see e.g. \cite{BlaCom,GhOs,Koh,FLan,Schoenb65,Schoenb66}) and
Sobolev's method \cite{Sobolev74,Sobolev06,SobVas}.  In different spaces,
based on these methods, Sard's problem was studied by many authors (see, for
example, \cite{IBab,BlaCom,CatCom,Com72a,Com72b,GhOs,HayMilShad11,Koh,FLan,MalOrl,Shad83,Shad02,ShadHay04,ShadHay11,ShadHayAzam12,SchSil,Sobolev06,Zag,ZhaSha} and references therein).

It should be noted that Sobolev's method is based on construction of
a discrete analogue to a linear differential operator.

In \cite{Sobolev74,Sobolev06,SobVas}, the minimization problem of the norm of the error functional by coefficients was reduced to the system of
difference equations of Wiener-Hopf type in the space $L_2^{(m)}$, where $L_2^{(m)}$ is the Sobolev space of functions with square integrable
generalized $m$th derivative. The existence and uniqueness of a solution of this system was proved
by Sobolev \cite{Sobolev74,SobVas}, who gave a description of some analytic algorithm for finding the coefficients of optimal
cubature formulas. For this, Sobolev defined and studied the discrete analogue $D^{(m)}_{hH}(h\beta)$ of the polyharmonic
operator $\Delta^m$. The problem of construction of the discrete operator $D^{(m)}_{hH}(h\beta)$ in $n$-dimensional
case was very hard and steel is one of open problems. In the one-dimensional case, the discrete analogue $D^{(m)}_h(h\beta)$
of the differential operator $\d^{2m}/\d x^{2m}$ was constructed by Z.Zh.~Zhamalov \cite{Zham78} and Kh.M.~Shadimetov \cite{Shad85}.
In the space $L_2^{(m)}$ using the discrete analogue $D^{(m)}_h(h\beta)$ of the differential operator $\d^{2m}/\d x^{2m}$ the following results were obtained:
optimal quadrature formulas of the form (\ref{(1)}) were constructed in \cite{Shad83} and positivity of their coefficients was studied in \cite{ShadHay11};
weight optimal quadrature formulas were obtained in the work \cite{Shad02};  optimal quadrature formulas of Euler-Maclaurin type were studied in \cite{ShadHayNur13} and \cite{ShadHayNur16}; the problem of construction of $D^m$-splines were considered in \cite{CabHayShad14};
optimal quadrature formulas for numerical integration of Fourier coefficients were constructed in the work \cite{BolHayShad17}.
We note that the above mentioned optimal quadrature formulas and splines, which are constructed in the space $L_2^{(m)}$, are exact for any algebraic polynomials of degree $m-1$.

Further, a discrete analogue of the differential operator $\frac{\d^{2m}}{\d x^{2m}}-\frac{\d^{2m-2}}{\d x^{2m-2}}$ was constructed
in \cite{ShadHay04a} and properties of this discrete analogue were studied in \cite{ShadHay04b}.
This discrete analogue was applied for finding the coefficients of the optimal quadrature formulas
(see \cite{ShadHay04,ShadHay14}), for construction of interpolation splines minimizing the semi-norm (see \cite{ShadHay13}) and
for obtaining explicit coefficients of optimal quadrature formulas for Fourier coefficients \cite{BolHayMilShad17} in the Hilbert space
$W_2^{(m,m-1)}$. It should be noted that the obtained optimal quadrature formulas
and splines in the space $W_2^{(m,m-1)}$ are exact for any algebraic polynomial of degree $m-2$ and for exponential function $\E^{-x}$.

Next, in the Hilbert space $K_2(P_m)$ optimal quadrature formulas \cite{HayMilShad11,HayMilShad15} and interpolation splines minimizing the semi-norm
\cite{Hay18,HayMilShad14} were obtained using the discrete analogue of the differential operator $\frac{\d^{2m}}{\d x^{2m}}+2\omega^2\frac{\d^{2m-2}}{\d x^{2m-2}}+\omega^4\frac{\d^{2m-4}}{\d x^{2m-4}}$ for $m\geq 2$ \cite{Hay14}. Note that the obtained formulas are exact for algebraic polynomials
of degree $m-3$ and for trigonometric functions $\sin\omega x$ and $\cos\omega x$.

We note that optimal quadrature formulas and interpolation splines are exact to the base functions of a kernel of the semi-norm given
in a Hilbert space. The above mentioned optimal quadrature formulas and splines constructed in the spaces $L_2^{(m)}$, $W_2^{(m,m-1)}$ and $K_2(P_m)$
are exact to an algebraic polynomial of some degree which is the base function of the kernel in the corresponding semi-norm.

The aim of this work is to study Sard's problem on construction of optimal quadrature formulas of the form (\ref{(1)})
in the space $W_2^{(m,0)}$ by Sobolev's method. As a result we get the optimal quadrature formula which is exact to the base functions of the
kernel of the norm (\ref{norm}). Here the base functions consist on exponential-trigonometric functions.
It may be noted that we obtain the optimal quadrature formula which is not exact for any algebraic polynomial.

The rest of the paper is organized as follows:  in Sect. 2 the extremal function, which corresponds
to the error functional $\ell$, is found and, in Sect. 3, with the help of this extremal function
the norm of the error functional is calculated, i.e., the first part of Sard's problem is solved; in Sect. 4
the system of linear equations for coefficients of the optimal quadrature
formulas in the space $W_2^{(m,0)}$ is obtained, moreover the existence and uniqueness of the
solution of this system are discussed; in Sect. 5 the discrete analogue $D_m[\beta]$ of the differential operator $\frac{\d^{2m}}{\d x^{2m}}-1$ is constructed;
in Sect. 6 Sobolev's method of construction of optimal quadrature formulas of the form (\ref{(1)}) in the space $W_2^{(m,0)}$ is discussed.
Finally, in Sect. 7, for $m=1$ and $m=3$, the optimal quadrature formulas of the form (\ref{(1)}) are
given.

\section{The extremal function}

For finding the norm of the error functional (\ref{(2)}) of the quadrature formula (\ref{(1)}) we use the extremal function
of the error functional.

The function $\psi_{\ell}$ satisfying
the equation
\begin{equation}\label{(6)}
\left( {\ell ,\psi_\ell} \right) = \left\| \ell\right\|_{W_2^{(m,0)*}}
\left\| \psi_\ell\right\|_{W_2^{(m,0)}},
\end{equation}
is called \emph{the extremal function} for the functional $\ell$  \cite{Sobolev74,SobVas}.

Since $W_2^{(m,0)}$ is the Hilbert space then, by the Riesz theorem on general form of a linear continuous
functional on Hilbert spaces, for the error functional $\ell\in W_2^{(m,0)*}$ there exists a unique function $\psi_{\ell}\in W_2^{(m,0)}$ such that
for any $\varphi\in W_2^{(m,0)}$  the following equality is fulfilled
\begin{equation}\label{(7)}
\left( {\ell,\varphi} \right) =  \langle\psi_\ell,\varphi\rangle_m.
\end{equation}
Moreover, $\left\| \ell\right\|_{W_2^{(m,0)*}} = \left\| \psi_\ell\right\|_{W_2^{(m,0)}}$. Here
$\langle\psi_\ell,\varphi\rangle_m$ is the inner product of two functions $\psi_{\ell}$ and $\varphi$ defined by equality (\ref{inner-product}) in the space $W_2^{(m,0)}$. In particular, from (\ref{(7)}) when $\varphi=\psi_{\ell}$ we have
$$
\left( {\ell,\psi_\ell} \right) = \langle\psi_\ell,\psi_\ell\rangle_m =
\left\| \psi_\ell \right\|_{W_2^{(m,0)}}^2  = \left\| {\ell} \right\|_{W_2^{(m,0)*}}^2.
$$
It is clear that the solution $\psi_{\ell}$ of equation (\ref{(7)}) is the extremal function.
Thus, in order to calculate the norm of the error functional $\ell$, first, we should find the extremal function $\psi_{\ell}$ from equation (\ref{(7)})
and then to calculate the square of the norm of the error functional $\ell$ as follows
\begin{equation}\label{(8)}
\left\| {\ell} \right\|_{W_2^{(m,0)*}}^2=(\ell,\psi_\ell).
\end{equation}

In the present paper we study Sard's problem in the space $W_2^{(m,0)}$ for odd natural numbers $m$.

Further, we suppose that $m$ is an odd natural number.

Integrating by parts the right hand side of (\ref{(7)})
we get
\begin{eqnarray*}
(\ell,\varphi)&=&(-1)^m \int\limits_0^1 \left( \psi_\ell ^{(2m)} (x) - \psi_\ell (x) \right)\varphi (x)\d x\\
&&+ \left. {\sum\limits_{s = 0}^{m - 1} {( - 1)^{s} } \left( {\psi_\ell
^{(m + s)} (x) + \psi_\ell ^{(s)} (x) } \right)\varphi ^{(m-s- 1)}(x)} \right|_0^1.
\end{eqnarray*}
From here taking into account uniqueness of the function $\psi_{\ell}$ we have the following equation
\begin{equation}\label{(9)}
\psi_\ell^{(2m)}(x)-\psi_\ell(x)=(-1)^m\ell(x)
\end{equation}
with the boundary conditions
\begin{equation}\label{(10)}
\left.\left( \psi_\ell^{(m+s)}(x) +
\psi_\ell^{(s)}(x)\right)\right|_{x=0}^{x=1}=0, \ s=0,1,...,m-1.
\end{equation}

We note that in the work \cite{Boltaev11} for the solution of the boundary value problem (\ref{(9)})-(\ref{(10)})
the following result was obtained.

\begin{theorem}\label{Thm1}
A solution of equation (\ref{(9)}) with the boundary conditions (\ref{(10)}) is the
extremal function $\psi_{\ell}$ of the error functional $\ell$ and has the form
\begin{equation}\label{(12)}
\psi_\ell(x) = (-1)^m \ell(x)*G_m(x)+Y(x),
\end{equation}
where
\begin{equation}\label{(13)}
G_m(x)=\frac{\mathrm{sgn}x}{2m} \left(\sinh(x)+\sum_{n=1}^{m-1}\E^{x  \cos\left(\frac{\pi n}{m}\right)}
\cos\left(x\sin\left(\frac{\pi n}{m}\right)+\frac{\pi n}{m}\right)\right)
\end{equation}
and
\begin{eqnarray}
Y(x)&=&d_0\E^{-x}+\sum\limits_{k=1}^{\frac{m-1}{2}}\E^{x\cos\frac{(2k-1)\pi }{m}}\Bigg(d_{1,k}\cos\left(x\sin
\left(\frac{(2k-1)\pi}{m}\right)\right)\nonumber \\
&&\qquad\qquad+ d_{2,k}\sin\left(x\sin\left(\frac{(2k-1)\pi}{m}\right)\right)\Bigg),\label{(14)}
\end{eqnarray}
$m$ is an odd natural number, $d_0,\ d_{1,k},$ and $d_{2,k}$ are constants.
\end{theorem}

Since the error functional (\ref{(2)}) is defined on the space $W_2^{(m,0)}$ it is necessary to impose
the following conditions
\begin{eqnarray}
&&(\ell,\E^{-x}  ) = 0,\label{(15)}  \\
&&\left({\ell,\E^{-x\cos\frac{2\pi k}{m}} \cos {\left(x\sin\frac{2\pi k}{m}\right)}}
\right) = 0,  \ k=1,2,...,\frac{m-1}{2},\label{(16)}\\
&&\left({\ell,\E^{-x\cos\frac{2\pi k}{m}} \sin {\left(x\sin\frac{2\pi k}{m}\right)}}
\right) = 0,\  k=1,2,...,\frac{m-1}{2},\label{(17)}
\end{eqnarray}
which mean that the quadrature formula (\ref{(1)}) is exact for linear combinations of functions
$$
\E^{-x},\ \E^{-x\cos\frac{2\pi k}{m}} \cos \left(x\sin\frac{2\pi k}{m}\right)\mbox{ and }
\E^{-x\cos\frac{2\pi k}{m}} \sin \left(x\sin\frac{2\pi k}{m}\right)
$$
for $k=1,2,...,\frac{m-1}{2}$.

It should be pointed out that to ensure the solvability of the system (\ref{(15)})-(\ref{(17)}) with respect
to the coefficients $C_{\beta}$ ($\beta=0,1,...,N$), the condition $N +1\geq m$ has to be imposed.

Further, we also use the following equivalent forms of conditions (\ref{(16)})
and (\ref{(17)}) which are obtained by changing of variable $k=\frac{m+1}{2}-j$ with $k=1,2,...,\frac{m-1}{2}$
\begin{eqnarray}
&&\left({\ell,\E^{x\cos\frac{(2j-1)\pi}{m}} \cos {\left(x\sin\frac{(2j-1)\pi}{m}\right)}}
\right) = 0,  \ j=1,2,...,\frac{m-1}{2},\label{(18)}\\
&&\left({\ell,\E^{x\cos\frac{(2j-1)\pi}{m}} \sin {\left(x\sin\frac{(2j-1)\pi}{m}\right)}}
\right) = 0,\  j=1,2,...,\frac{m-1}{2}.\label{(19)}
\end{eqnarray}
Now, in the next section, using the extremal function (\ref{(12)}), we calculate the norm of the error functional $\ell$.

\section{The norm of the error functional $\ell$}

As was said above, in order to calculate the square of the norm of the error functional (\ref{(2)}) it is sufficient to calculate the value $(\ell,\psi_\ell)$
of the error functional $\ell$ at the function $\psi_\ell$.
For this, firstly, using equalities (\ref{(15)}), (\ref{(18)}) and (\ref{(19)}), we get
$$
(\ell,Y(x))=0,
$$
where $Y(x)$ is the function defined by (\ref{(14)}).
Then, using (\ref{(12)}), taking the last equality into account and keeping in mind that $m$ is an odd natural number, we have
\begin{eqnarray}
\left\|\ell \right\|^2=\left(\ell,\psi_\ell\right)& =&
-\int\limits_{-\infty}^{\infty}\ell(x)\left[\ell(x)*G_m(x) + Y(x)\right]\d x\nonumber \\
&=&-\int\limits_{-\infty}^{\infty}\ell (x)\left[ \ell (x)*G_m (x)\right]\d x.\label{(20)}
\end{eqnarray}
Now, for the convolution in (\ref{(20)}), taking (\ref{(2)}) into account, we obtain
$$
\ell(x) * G_m (x)= \int\limits_{-\infty}^{\infty}\ell(y)G_m(x - y)dy=
\int\limits_0^1  G_m (x - y)dy - \sum\limits_{\beta  = 0}^N C_\beta  G_m (x - x_\beta).
$$
Then the square of the norm (\ref{(20)}) of the error functional $\ell$ takes the form
\begin{eqnarray}
\left\|\ell\right\|^2&=&\sum\limits_{\beta=0}^N
C_\beta  \int\limits_0^1 \left(G_m (x - x_\beta)+G_m (x_\beta - x)\right)\d x\nonumber\\
&&-\sum\limits_{\beta  = 0}^N \sum\limits_{\gamma  = 0}^N
C_\beta  C_\gamma G_m (x_\beta   - x_\gamma)-
\int\limits_0^1 \int\limits_0^1 G_m (x - y)\d x\d y.\label{(21)}
\end{eqnarray}
Since $G_m(x)$, defined by (\ref{(13)}), is the even function we have
$$
G_m(x_\beta-x)=G_m(x - x_\beta).
$$
Then, taking into account the last equality, from (\ref{(21)})
we get
\begin{eqnarray}
\left\|\ell \right\|^2 & =&2\sum\limits_{\beta  = 0}^N
C_\beta \int\limits_0^1 G_m (x - x_\beta)\d x \nonumber\\
&&-\sum\limits_{\beta  = 0}^N \sum\limits_{\gamma  = 0}^N
C_\beta  C_\gamma G_m (x_\beta   - x_\gamma)- \int\limits_0^1 \int\limits_0^1 G_m (x - y)\d x\d y.\label{(22)}
\end{eqnarray}
Thus, the first part of Sard's problem on construction of optimal quadrature formulas in the space $W_2^{(m,0)}$
is solved. Next we consider the second part of the problem.

\section{The discrete system of Wiener-Hopf type}

We consider the problem of minimization of the square of the norm (\ref{(22)}) for the error functional $\ell$ by coefficients $C_\beta$.
The error functional $\ell$ satisfies the conditions (\ref{(15)})-(\ref{(17)}).
The square of the norm (\ref{(22)}) of the error functional $\ell$ is a multi-variable quadratic function with respect to coefficients $C_\beta$ of the quadrature
formula (\ref{(1)}). For finding a point of conditional minimum of the function (\ref{(22)}) under the conditions
(\ref{(15)})-(\ref{(17)}) we apply the Lagrange method.

Denoting as $\mathbf{C }= (C_0 ,C_1,...,C_N )$ and  $\mathbf{d}=(d_0 ,\,d_{1,1},...,
d_{1,\frac{m-1}{2}},d_{2,1},...,d_{2,\frac{m-1}{2}})$ we consider the following function
\begin{eqnarray*}
\Psi(\mathbf{C}, \textbf{d})&=&\left\|\ell \right\|^2  - 2d_0 \left( {\ell,\E^{-x}} \right)-2\sum\limits_{k=1}^{\frac{m-1}{2}} \Bigg[ d_{1,k} \left(\ell,\E^{-x\cos\frac {2\pi k}{m}} \cos \left(x\sin\frac {2 \pi k}{m} \right) \right)\\
&&\qquad\qquad\qquad+ d_{2,k} \left(\ell,\E^{-x\cos\frac {2\pi k}{m}} \sin \left(x\sin\frac {2 \pi k}{m} \right) \right)\Bigg],
\end{eqnarray*}
where $d_0, d_{1,k}$ and $d_{2,k}$ ($k=1,2,...,\frac{m-1}{2}$) are Lagrange multipliers.

Equating to zero the partial derivatives of the function $\Psi(\mathbf{C}, \textbf{d})$ by coefficients $C_{\beta}$, $\beta=0,1,...,N$
and by $d_0, d_{1,k}$, $d_{2,k}$, $k=1,2,...,\frac{m-1}{2}$ we get the following difference system of $N+1+m$ linear equations with $N+1+m$ unknowns
\begin{eqnarray}
&&\sum_{\gamma=0}^N C_\gamma  G_m(x_\beta-x_\gamma)+P(x_\beta,d_0,d_{1,k},d_{2,k})= f_m(x_\beta),
\ \ \   \beta  = 0,1,\ldots,N,\label{(23)} \\
&&\sum_{\gamma  = 0}^N C_\gamma  \E^{-x_\gamma}  = 1 -\frac{1}{\E},\label{(24)}\\
&&\sum_{\gamma  = 0}^N C_\gamma \E^{-x_\gamma \cos \frac{2 \pi k}{m}}
\cos \left( x_\gamma \sin \frac{2 \pi k}{m}\right) = g_{1,k}, \ k=1,2,...,\frac{m-1}{2}, \label{(25)}\\
&&\sum_{\gamma  = 0}^N C_\gamma \E^{-x_\gamma \cos \frac{2 \pi k}{m}}
\sin \left( x_\gamma \sin \frac{2 \pi k}{m}\right) = g_{2,k},\ k=1,2,...,\frac{m-1}{2},\label{(26)}
\end{eqnarray}
where
\begin{eqnarray}
&&P(x_\beta,d_0,d_{1,k},d_{2,k})=d_0 \E^{-x_\beta} \nonumber \\
&&\qquad + \sum_{k=1}^{\frac{m-1}{2}} \E^{-x_\beta \cos \frac{2 \pi k}{m}}
\left( d_{1,k}\cos \left(x_\beta\sin \frac{2 \pi k}{m}\right)+d_{2,k}\sin \left(x_\beta\sin \frac{2 \pi k}{m}\right) \right),\label{(27)}\\
&&f_m(x_\beta) =\int_0^1 {G_m (x - x_\beta  )\d x}, \label{(28)}\\
&&g_{1,k}=\cos \frac {2 \pi k}{m}-\E^{-\cos \frac {2 \pi k}{m}}  \cos
\left( \sin \frac {2 \pi k}{m}+\frac {2 \pi k}{m}\right),\ k=1,2,...,\frac{m-1}{2},\label{(29)}\\
&&g_{2,k} = \sin \frac {2 \pi k}{m}-\E^{-\cos \frac {2 \pi k}{m}}  \sin
\left( \sin \frac {2 \pi k}{m}+\frac {2 \pi k}{m}\right),\ k=1,2,...,\frac{m-1}{2} \label{(30)}
\end{eqnarray}
and $G_m(x)$ is defined by (\ref{(13)}).

We note that the system (\ref{(23)})-(\ref{(26)}) is called \emph{the discrete system of Wiener-Hopf type }\cite{Sobolev74,SobVas}.

It is important to note that the existence and uniqueness of an optimal quadrature formula of the form (\ref{(1)})
in the sense of Sard in Hilbert spaces were studied in \cite{FLan}. Particularly,
we get that the difference system (\ref{(23)})-(\ref{(26)}) for any set of different nodes $x_{\beta}$, $\beta=0,1,...,N$ when $N+1\geq m$ has a unique solution and this solution gives minimum to $\|\ell\|^2$ defined by (\ref{(22)}) under the conditions (\ref{(15)})-(\ref{(17)}).
Existence and uniqueness of the solution of such a type of difference systems were also studied in \cite{HayMilShad11,ShadHay14,Sobolev74,SobVas}.

Further, we consider the case of equally spaced nodes. Suppose $x_\beta=h\beta$, $\beta=0,1,...,N$, $h=\frac{1}{N}$, $N=1,2,...$.

We assume that $C_\beta=0$ for $\beta<0$ and $\beta>N$. Then, using the convolution of two discrete argument functions
$\varphi(h\beta)$ and $\psi(h\beta)$ (see \cite{Sobolev74,SobVas})
$$
\varphi(h\beta)*\psi(h\beta)=\sum\limits_{\gamma=-\infty}^\infty  {\varphi(h\gamma) \cdot \psi(h\beta  -h \gamma}),
$$
we rewrite the system (\ref{(23)})-(\ref{(27)}) in the following convolution form
\begin{eqnarray}
&&C_\beta*G_m(h\beta)+P(h\beta,d_0,d_{1,k},d_{2,k})=f_m(h\beta),\ \ \   \beta=0,1,...,N, \label{(31)}\\
&&\sum_{\gamma  = 0}^N C_\gamma  \E^{-h\gamma}  = 1 -\frac{1}{\E},\label{(32)}\\
&&\sum_{\gamma  = 0}^N C_\gamma \E^{-h\gamma \cos \frac{2 \pi k}{m}}\cos \left( h\gamma \sin \frac{2 \pi k}{m}\right) = g_{1,k},
\ k=1,2,...,\frac{m-1}{2},\label{(33)}\\
&&\sum_{\gamma  = 0}^N C_\gamma \E^{-h\gamma \cos \frac{2 \pi k}{m}}\sin \left( h\gamma \sin \frac{2 \pi k}{m}\right) = g_{2,k},
\ k=1,2,...,\frac{m-1}{2},\label{(34)}
\end{eqnarray}
where $G_m(h\beta)$, $P(h\beta,d_0,d_{1,k},d_{2,k})$, $f_m(h\beta)$, $g_{1,k}$ and $g_{2,k}$ are defined by (\ref{(13)}), (\ref{(27)})-(\ref{(30)}),
respectively. In the system (\ref{(31)})-(\ref{(34)}) there are $N+m+1$ unknowns $C_\beta,$ $\beta=0,1,...,N$, $d_0,d_{1,k}$ and $d_{2,k}$,
$k=1,2,...,\frac{m-1}{2}$ and $N+m+1$ linear equations.

In order to solve the system (\ref{(31)})-(\ref{(34)}) by Sobolev's method we need discrete analogue of the differential operator
$\frac{\d^{2m}}{\d x^{2m}}-1$. The next section is devoted to construction of this discrete analogue.

\section{A discrete analogue of the differential operator $\frac{\d^{2m}}{\d x^{2m}}-1$}

In the present section for an odd natural number $m$ we construct the function $D_m(h\beta)$ of discrete argument
which satisfies the equation
\begin{equation}\label{(35)}
D_m(h\beta)*G_m(h\beta) = \delta_{\d}(h\beta),
\end{equation}
where
\begin{equation}\label{(36)}
G_m(h\beta)={{{\rm sgn}(h\beta)} \over 2m}
 \left( \sinh(h\beta)+\sum\limits_{n=1}^{m-1}{\E^{h\beta\cos{\frac{\pi n}{m}}}\cos{\left(h\beta\sin{\frac{\pi n}{m}+\frac{\pi n}{m}}\right)}}\right), \end{equation}
$\delta_{\d}(h\beta)$ is the discrete delta-function, i.e., $\delta_{\d}(h\beta)=\left\{
\begin{array}{ll}
1, & \beta=0,\\
0, &\beta\neq 0.
\end{array}
\right. $

It should be noted that a method of construction of the discrete argument function $D_m(h\beta)$
similarly as the method of construction of discrete analogues of the differential operators
$\frac{\d^{2m}}{\d x^{2m}}$, $\frac{\d^{2m}}{\d x^{2m}}-\frac{\d^{2m-2}}{\d x^{2m-2}}$ and
 $\frac{\d^{2m}}{\d x^{2m}}+2\omega^2\frac{\d^{2m-2}}{\d x^{2m-2}}+\omega^4\frac{\d^{2m-4}}{\d x^{2m-4}}$
in the works \cite{Shad85,ShadHay04a} and \cite{Hay14}.

We note that equation (\ref{(35)}) is a discrete analogue of the following equation
\begin{equation}\label{(37)}
\left(\frac{\d^{2m}}{\d x^{2m}}-1\right)G_m(x) = \delta(x),
\end{equation}
where $G_m(x)$ is defined by (\ref{(13)}), $\delta$ is Dirac's delta-function.

The results of this section are the following.

\begin{theorem}\label{Thm2}
The discrete analogue of the differential operator
${{\d^{2m} } \over {\d x^{2m}}} - 1$ satisfying the equation (\ref{(35)}),
when $m$ is an odd natural number, has the form
\begin{equation}\label{(38)}
D_m(h\beta)={\frac{m}{K}}\left\{
\begin{array}{ll}
\sum\limits_{n=1}^{m-1}A_n \lambda_n^{|\beta|-1},& |\beta|\geq 2,\\
1+\sum\limits_{n=1}^{m-1}A_n, & |\beta|=1,\\
\ M_1-\frac{K_1}{K}+\sum\limits_{n=1}^{m-1}\frac{A_n}{\lambda_n},& \beta=0,
\end{array}
\right.
\end{equation}
where
\begin{eqnarray*}
A_n&=&\frac{(\lambda_n^2-2\lambda_n \cosh(h)+1)B_{2m-2}(\lambda_n)}{\lambda_n {\cal P}_{2m-2}'(\lambda_n)},\ n=1,2,...,m-1,\\
K&=&\sum\limits_{k=1}^{\frac{m-1} {2}} a_{1,k}+\sinh(h),\\
K_1 &=&\sum\limits_{k=1}^{\frac{m-1} {2}}\left(b_{1,k}\sinh(h)+a_{2,k}+a_{1,k}\left(\sum\limits_{j=1,j\neq k}^{\frac{m-1} {2}} b_{1,j}-2\cosh(h)\right)\right),\\
M_1& =&\sum\limits_{k=1}^{\frac{m-1}{2}} b_{1,k}-2\cosh(h),\\
B_{2m-2}(\lambda)&=&\prod\limits_{k=1}^{\frac{m-1}{2}}{\left( {{\lambda }^{4}}+b_{1,k}{{\lambda }^{3}}
+b_{2,k}{{\lambda }^{2}}+b_{1,k}\lambda+1 \right)},\\
{\cal{P}}_{2m-2}(\lambda)&=&\left(\sinh(h)+(\lambda^2-2\lambda \cosh(h)+1)\sum\limits_{j=1}^{\frac{m-1}{2}}\frac{a_{1,j}\lambda^2+a_{2,j}\lambda+a_{1,j}}
{\lambda^4+b_{1,j}\lambda^3+b_{2,k}\lambda^2+b_{1,k}\lambda+1}\right)B_{2m-2}(\lambda)
\end{eqnarray*}
and here for $k=1,2,...,\frac{m-1}{2}$
\begin{eqnarray*}
a_{1,k} &=& 2\cdot\left[\cos{\frac{\pi k}{m}}\cos{\left(h\sin{\frac{\pi k}{m}}\right)}\sinh{\left(h\cos{\frac{\pi k}{m}}\right)}\right.\\
  & &\qquad\qquad\left.-\sin{\frac{\pi k}{m}}\sin{\left(h\sin{\frac{\pi k}{m}}\right)}\cosh{\left(h\cos{\frac{\pi k}{m}}\right)}\right],\\
a_{2,k} &=& -2\cdot\left[\cos{\frac{\pi k}{m}}\sinh{\left(2h\cos{\frac{\pi k}{m}}\right)}-\sin{\frac{\pi k}{m}}\cosh{\left(2h\cos{\frac{\pi k}{m}}\right)}\right],\\
b_{1,k} &=& -4\cdot\cos{\left(h\sin{\frac{\pi k}{m}}\right)}\cosh{\left(h\cos{\frac{\pi k}{m}}\right)},\\
b_{2,k} &=& 2\cdot\left[1+\cos{\left(2h\sin{\frac{\pi k}{m}}\right)}+\cosh{\left(2h\cos{\frac{\pi k}{m}}\right)}\right],\\
\end{eqnarray*}
$\lambda_n$, $n=1,2,...,m-1$ are roots of the polynomial $\mathcal{P}_{2m-2}(\lambda)$ with $|\lambda_n|<1$, $h$ is a small positive parameter.
\end{theorem}

\begin{theorem}\label{Thm3}
The discrete analogue $D_m(h\beta)$ of the differential operator
${{\d^{2m}}\over{\d x^{2m}}} - 1$ satisfies the following equalities\\
\begin{tabular}{ll}
1) $D_m\left(h\beta\right)*{{\E}^{ h\beta  }}=0,$ &\\
2) $D_m\left(h \beta  \right)*{{\E}^{-h \beta}}=0,$&\\
3) $D_m\left(h\beta  \right)*{{\E}^{h \beta \cos \left( \frac{2\pi k}{m} \right)}} \cos \left( h\beta \sin \frac{2\pi k}{m} \right)=0,$ & \\
4) $D_m\left(h \beta  \right)*{{\E}^{-h \beta \cos \left( \frac{2\pi k}{m} \right)}}  \cos \left( h\beta\sin \frac{2\pi k}{m} \right)=0,$ & \\
5) $D_m\left(h \beta  \right)*{{\E}^{h\beta \cos \left( \frac{2\pi k}{m} \right)}}  \sin \left( h\beta \sin \frac{2\pi k}{m} \right)=0,$& \\
6) $D_m\left(h\beta  \right)*{{\E}^{-h \beta \cos \left( \frac{2\pi k}{m} \right)}}  \sin \left( h \beta \sin \frac{2\pi k}{m} \right)=0,$& \\
\end{tabular}\\
where  $k=1,2,...,\frac{m-1}{2},$ $m$ is an odd natural number.
\end{theorem}

Here we give the proof of Theorem \ref{Thm2}. Theorem \ref{Thm3} can be proved by direct calculation of convolutions
on the left hand sides of equalities 1) - 6).

\medskip

For proving Theorems \ref{Thm2} we need the following well-known formulas from the theory of generalized functions and
Fourier transforms (see, for instance, \cite{Sobolev74,Vlad})
\begin{equation}\label{(39)}
F\left[ \phi \left(x\right) \right]=\int\limits_{-\infty }^{\infty }{\phi \left( x \right){{e}^{2\pi \i px}}\d x,\,\,\,\,}
{{F}^{-1}}\left[ \phi \left( p \right)\right]=\int\limits_{-\infty }^{\infty }{\phi \left( p \right){{e}^{-2\pi \i px}}\d p},
\end{equation}
\begin{equation}\label{(40)}
F\left[\phi*\psi \right]=F\left[ \phi  \right]\cdot F\left[\psi\right],
\end{equation}
\begin{equation}\label{(41)}
F\left[ \phi \cdot \psi  \right]=F\left[ \phi  \right]*F\left[ \psi  \right],
\end{equation}
\begin{equation}\label{(42)}
F\left[\delta\left(x\right)\right]=1,\ \ F\left[\delta^{(\alpha)}(x)\right]={{\left(-2\pi \i p\right)}^{\alpha }},
\end{equation}
\begin{equation}\label{(43)}
\delta \left(hx \right)={{h}^{-1}}\delta \left(x\right),
\end{equation}
\begin{equation}\label{(44)}
\delta \left(x-a\right) f\left(x\right)=\delta \left(x-a\right) f\left( a \right),
\end{equation}
\begin{equation}\label{(45)}
\delta^{(\alpha)}(x)*f(x)=f^{(\alpha)}(x),
\end{equation}
\begin{equation}\label{(46)}
{{\phi }_{0}}\left( x \right)=\sum\limits_{\beta =-\infty }^{\infty }{\delta \left( x-\beta  \right)},\ \ \
\sum\limits_{\beta=-\infty}^{\infty}{{{e}^{2\pi \i x\beta}}}=\sum\limits_{\beta=-\infty}^{\infty}
{\delta \left( x-\beta  \right)}.
\end{equation}
\medskip

\textbf{Proof of Theorem \ref{Thm2}}. Based on the theory of generalized functions and Fourier transforms it is convenient
to find a harrow-shaped function instead of the discrete function $D_m(h\beta)$ (see \cite{Sobolev74,SobVas}) as follows
$$
\stackrel{\abd}{D}_m\!\!(x) = \sum\limits_{\beta  = - \infty }^\infty  {D_m(h\beta)\delta } (x - h\beta ).
$$
In the class of harrow-shaped functions the equation (\ref{(35)})
is reduced to the following equation
\begin{equation}\label{(47)}
\stackrel{\abd}{D}_m\!\!(x)*\stackrel{\abd}{G}_m\!\!(x)= \delta (x),
\end{equation}
where $\stackrel{\abd}{G}_m\!\!(x) = \sum\limits_{\beta  = - \infty }^\infty  {G_m(h\beta) \delta } (x - h\beta )$ is the harrow-shaped
function corresponding to the discrete function $G_m(h\beta)$.

Applying the Fourier transform to both sides of equation (\ref{(47)}), taking (\ref{(40)}) and (\ref{(42)}) into account,
we get
\begin{equation}\label{(48)}
F[\stackrel{\abd}{D}_m\!\!(x)] = {1 \over {F[\stackrel{\abd}{G}_m\!\!(x)]}}.
\end{equation}

First we calculate the Fourier transform $F[\stackrel{\abd}{G}_m\!\!(x)]$.
For this, using equalities (\ref{(44)}), (\ref{(46)}), (\ref{(43)}) and
(\ref{(41)}), we have
\begin{equation}\label{(49)}
F[\stackrel{\abd}{G}_m\!\!(x)] =  F[G_m(x)]*\phi_0(hp).
\end{equation}
Using equalities (\ref{(37)}) and (\ref{(45)}), keeping in mind (\ref{(40)}) and (\ref{(42)}),
we get
\begin{equation}\label{(50)}
F[G_m(x)] = {1 \over {(2\pi \i p)^{2m}-1}}.
\end{equation}
We denote roots of the equation $p^{2m}-1=0$ as
$$
p_1 = 1,\ p_2 = -1,\ p_{1,n}=\cos\left(\frac{\pi n}{m}\right)+\i\sin\left(\frac{\pi n}{m}\right),
\mbox{ and }p_{2,n}  = \cos\left(\frac{\pi n}{m}\right)-\i\sin\left(\frac{\pi n}{m}\right),
$$ where  $n=1,2,...,m-1$.

Then calculation of the convolution in the right hand side of (\ref{(49)}), by taking (\ref{(48)}) and (\ref{(50)}) into account, gives
\begin{eqnarray}
F[\stackrel{\abd}{D}_m](p) &=& -{{(2\pi )^{2m} } \over {h^{2m-1} }}
\Bigg[ \sum\limits_{\beta=-\infty}^{\infty}
{1 \over {[\beta  - h(p + {{\i p_1 } \over {2\pi }})] [\beta - h(p + {{\i p_2 } \over {2\pi }})]}}\nonumber\\
&&\times \prod\limits_{n=1}^{m-1}\frac{1}{[\beta  - h(p + {{\i p_{1,n} }
 \over {2\pi }})] [\beta  - h(p + {{\i p_{2,n} } \over {2\pi }})]} \Bigg]^{ - 1}.\label{(51)}
\end{eqnarray}
Suppose that the Fourier series of the function $F[\stackrel{\abd}{D}_m](p)$ has the form
\begin{equation}\label{(52)}
F[\stackrel{\abd}{D}_m](p) = \sum\limits_{\beta= - \infty }^\infty
 {\widehat D_m(h\beta)\E^{2\pi \i ph\beta } },
\end{equation}
where $\widehat D_m(h\beta)$ is the Fourier coefficients of the function $F[\stackrel{\abd}{D}_m](p)$, i.e.,
\begin{equation}\label{(53)}
\widehat D_m(h\beta) = \int\limits_0^{h^{ - 1} } {F[\stackrel{\abd}{D}_m](p)\ \E^{ - 2\pi \i ph\beta }
 } \d p.
\end{equation}
Applying the inverse Fourier transform to both sides of equality (\ref{(52)}) we get the harrow-shaped function
$$
\stackrel{\abd}{D}_m\!\!(x) = \sum\limits_{\beta  = - \infty }^\infty  {\widehat D_m(h\beta) \delta } (x - h\beta ).
$$
Hence, taking into account the definition of harrow-shaped functions, we conclude that the discrete function
$\widehat D_m(h\beta)$ is the discrete argument function $D_m(h\beta)$ which we are looking for.
This means that we can find $D_m(h\beta)$ by expanding the right hand side of (\ref{(51)}) into the Fourier series.

For calculation the series in (\ref{(52)}) we use the following well-known formula from the residual theory
(see, for instance, \cite[p. 296]{MaqSalSir76})
\begin{equation}\label{(54)}
\sum\limits_{\beta=-\infty }^\infty{f(\beta )=-\sum\limits_{z_1 ,z_2,...,z_t}{\rm res} (\pi \cot(\pi z) f(z))},
\end{equation}
where $z_1 ,z_2,...,z_t$ are the poles of the function $f(z)$.

We denote
$$
g(z) = {1 \over {[z - h(p + {{\i p_1} \over {2\pi }})]   [z - h(p + {{\i p_2}
\over {2\pi }})] \prod\limits_{n=1}^{m-1}[z - h(p + {{\i p_{1,n}} \over {2\pi }})]}[z - h(p + {{\i p_{2,n}} \over {2\pi }})]},
$$
where $z_1  = h(p + {{\i p_1} \over {2\pi }}),$  $z_2 = h(p + {{\i p_2} \over {2\pi }}),$  $z_{1,n}  = h(p + {{\i p_{1,n}} \over {2\pi }}),$
$z_{2,n}  = h(p + {{\i p_{2,n}} \over {2\pi i}})$, $n=1,2,...,m-1$ are the poles of order 1 of the function $g(z)$.
Then using the formula (\ref{(54)}) from (\ref{(51)}) we can write
\begin{equation}\label{(55)}
F[\stackrel{\abd}{D}_m](p)  ={{(2\pi)^{2m}}\over {h^{2m-1} }}
\left[\sum\limits_{z_1,z_2,z_{1,1},...,z_{1,m-1},z_{2,1},...,z_{2,m-1}}
{\rm{res}}(\pi \cot(\pi z)  g(z)) \right]^{ - 1}.
\end{equation}
By direct calculation we have the following results
\begin{eqnarray*}
\mathop{\rm{res}}\limits_{z = z_1}(\pi \cot(\pi z) g(z))
&=& \frac{{{\left( 2\pi  \right)}^{2m}}\cot\left(\pi hp+\frac{{{p}_{1}}h\i}{2} \right)}{4m {{\left( h\i \right)}^{2m-1}}},\\
\mathop {\mbox{res}}\limits_{z = z_2 } (\pi \cot(\pi z) g(z))
&=& -\frac{{{\left( 2\pi  \right)}^{2m}}\cot\left( \pi hp+\frac{{{p}_{2}}h\i}{2} \right)}{4m {{\left( h\i \right)}^{2m-1}}},\\
\mathop {\mbox{res}}\limits_{z = z_{1,n} } (\pi \cot(\pi z) g(z))
&=& \frac{{{\left( 2\pi  \right)}^{2m}}\cot\left( \pi hp+\frac{{{p}_{1,n}}h\i}{2} \right)}{4m {{\left( h\i \right)}^{2m-1}}\left( \cos \left( \frac{\pi }{m} \right)-\i\sin \left( \frac{\pi }{m} \right)\right)},\ n=1,2,...,m-1,\\
\mathop {\mbox{res}}\limits_{z = z_{2,n} } (\pi \cot(\pi z) g(z))
&=&\frac{{{\left( 2\pi  \right)}^{2m}}\cot\left( \pi hp+\frac{{{p}_{2,n}}h\i}{2} \right)}{4m {{\left( h\i \right)}^{2m-1}}\left( \cos \left( \frac{\pi }{m} \right)+\i\sin \left( \frac{\pi}{m} \right) \right)},\ n=1,2,...,m-1.
\end{eqnarray*}

Denoting by $\lambda ={{e}^{2\pi \i ph}}$, using the last $2m-2$ equalities and taking into account the following formulas
$$
\cos z=\frac{{{e}^{z\i}}+{{e}^{-z\i}}}{2},\,\,\sin z=\frac{{{e}^{z\i}}-{{e}^{-z\i}}}{2\i},\,\,\cosh z=\frac{{{e}^{z}}+{{e}^{-z}}}{2},\mbox{ and }\sinh z=\frac{{{e}^{z}}-{{e}^{-z}}}{2},
$$
after some calculations from (\ref{(55)}) we get
\begin{equation}\label{(56)}
F[\stackrel{\abd}{D}_m](p)=\frac{m}{K}\frac{\left( {{\lambda }^{2}}-2\lambda \cosh(h) +1 \right)B_{2m-2}(\lambda)}{\lambda \mathcal{P}_{2m-2}(\lambda)},
\end{equation}
where $K$, $B_{2m-2}(\lambda)$ and $\mathcal{P}_{2m-2}(\lambda)$ are defined in the statement of Theorem \ref{Thm2}.

For finding the explicit form of the discrete function $D_m[\beta]$ the right hand side of (\ref{(56)})
we expand into the sum of partial fractions as follows
\begin{eqnarray}
&& \frac{m}{K} \frac{\left( {{\lambda }^{2}}-2\lambda \cosh\left( h \right)+1 \right)B_{2m-2}(\lambda)}{\lambda \left( \lambda -{{\lambda }_{1}} \right)\left( \lambda -{{\lambda }_{2}} \right)\cdot ...\cdot \left( \lambda -{{\lambda }_{2m-2}} \right)}\nonumber\\
&&\qquad \left.=\frac{m}{K}\left( \lambda -\frac{{K}_{1}}{K}+{{M}_{1}}+\frac{{{A}_{0}}}{\lambda }\right.
+\frac{{{A}_{1}}}{\lambda -{{\lambda }_{1}}}+\frac{{{A}_{2}}}{\lambda -{{\lambda }_{2}}}+...+\frac{{{A}_{2m-2}}}{\lambda -{{\lambda }_{2m-2}}} \right),\label{(58)}
\end{eqnarray}
where $K$, $K_1$, $M_1$ and $B_{2m-2}(\lambda)$ are given in the statement of the theorem,  ${{A}_{0}},{{A}_{1}},...,{{A}_{2m-2}}$ are unknowns,
$\lambda_1,\lambda_2,...,\lambda_{2m-2}$ are roots of the polynomial $\mathcal{P}_{2m-2}(\lambda)$ such that
\begin{equation}\label{(57)}
\lambda_n\cdot \lambda_{2m-1-n}=1\mbox{ and }|\lambda_n|<1,\ |\lambda_{2m-1-n}|>1\mbox{ for }n=1,2,...,m-1.
\end{equation}
For finding unknowns ${{A}_{0}},{{A}_{1}},...,{{A}_{2m-2}}$ we multiply both sides of
equality (\ref{(58)}) by expression $\lambda(\lambda-\lambda_1)(\lambda-\lambda_2)\cdot ...\cdot (\lambda-\lambda_{2m-2})$
and we put $\lambda=0,\ \lambda=\lambda_1,...,\lambda=\lambda_{2m-2}$. Then we get the following results
\begin{eqnarray}
{{A}_{0}}&=&1,\label{(59)}\\
A_n&=&\frac{\left(\lambda _n^2-2\lambda_n\cosh(h)+1\right)B_{2m-2}(\lambda_n)}{\lambda_n\mathcal{P}_{2m-2}'(\lambda_n)}, \ \ n=1,2,...,2m-2. \label{(60)}
\end{eqnarray}
From (\ref{(60)}), taking (\ref{(57)}) into account, we have
\begin{equation}\label{(61)}
A_{2m-1-n}=-\frac{1}{\lambda _n^{2}}A_n,\ \ n=1,2,...,m-1.
\end{equation}
Finally, using (\ref{(57)})-(\ref{(61)}) from (\ref{(56)}) we get
\begin{eqnarray*}
F[\stackrel{\abd}{D}_m](p)&=&\frac{m}{K} \bigg[ \lambda -\frac{{K}_{1}}{K}+{{M}_{1}}+\frac{1}{\lambda } +\sum\limits_{\beta=0}^{\infty}
\sum\limits_{n=1}^{m-1}
\left(\frac{A_n}{\lambda}\left(\frac{\lambda_n}{\lambda}\right)^{\beta}+A_n\lambda_n(\lambda_n\lambda)^{\beta}\right)\bigg]\\
&=&\sum\limits_{\beta =-\infty }^{\infty }D_m(h\beta)\lambda^\beta.
\end{eqnarray*}
Hence, keeping in mind that $\lambda=\E^{2\pi\i ph}$, we obtain the explicit form (\ref{(38)}) of the discrete function $D_m(h\beta)$.
Theorem \ref{Thm2} is proved. \hfill $\Box$

\medskip

\begin{remark}
It is easy to see from (\ref{(38)}) that the function $D_m(h\beta)$ is an even function and
it decreases exponentially as $|\beta|\to \infty$.
\end{remark}

\section{Solution of the discrete Wiener-Hopf system}

In this section we give an algorithm for finding the exact solution of the system (\ref{(31)})-(\ref{(34)}) by using the discrete analogue $D_m(h\beta)$ of the differential operator $\frac{\d^{2m}}{\d x^{2m}}-1$ obtained in the previous section.

We introduce the following functions
\begin{eqnarray}
v_m(h\beta) &=& C_\beta*G_m(h\beta),\label{(62)}\\
u_m(h\beta) &=& v_m(h\beta) + P(h\beta, d_0, d_{1,k}, d_{2,k}).\label{(63)}
\end{eqnarray}
Then, taking (\ref{(35)}) into account, for optimal coefficients we have
\begin{equation}\label{(64)}
C_\beta = D_m(h\beta)*u_m(h\beta).
\end{equation}
Thus, if we find the function $u_m(h\beta)$ then optimal coefficients can be defined by the formula (\ref{(64)}).
For calculation of the convolution in (\ref{(64)}) it is required to find the function $u_m(h\beta)$ at all integer values of $\beta$.

It is clear from (\ref{(31)}) that $u_m(h\beta)=f_m(h\beta)$ for $\beta=0,1,...,N$.
Now we need to find the function $u_m(h\beta)$ for $\beta<0$ and $\beta>N$. Using the formula (\ref{(13)}),
we calculate the convolution $v_m(h\beta) = C_\beta*G_m(h\beta)$  for $\beta<0$ and $\beta>N$.\\
For $\beta<0$ we have
\begin{eqnarray}\label{(65)}
v_m(h\beta) &=&-\frac{1}{2m}\sum_{\gamma = 0 }^N C_\gamma \Bigg[\frac{\E^{h\beta  - h\gamma } - \E^{ - h\beta  + h\gamma }}{2}\nonumber\\
&&+\sum_{n = 1 }^{m-1} \E^{(h\beta-h\gamma) \cos \frac{\pi n}{m}}
\cos \left((h\beta-h\gamma) \sin \frac{\pi n}{m}+\frac{\pi n}{m}\right) \Bigg].
\end{eqnarray}
Hence splitting up the internal sum  into two parts by even and odd values of $n$, respectively, we get
\begin{equation}\label{(66)}
\sum_{n = 1 }^{m-1} \E^{(h\beta-h\gamma) \cos \frac{\pi n}{m}}
\cos \left((h\beta-h\gamma) \sin \frac{\pi n}{m}+\frac{\pi n}{m}\right) = S_{1}+S_{2},
\end{equation}
where
\begin{eqnarray*}
S_{1}& =& \sum_{k = 1 }^{\frac{m-1}{2}} \E^{(h\beta-h\gamma) \cos \frac{2\pi k}{m}}
\cos \left((h\beta-h\gamma) \sin \frac{2\pi k}{m}+\frac{2\pi k}{m}\right), \\
S_{2} &=& \sum_{k = 1 }^{\frac{m-1}{2}} \E^{(h\beta-h\gamma) \cos \frac{(2k-1)\pi}{m}}
\cos \left((h\beta-h\gamma) \sin \frac{(2k-1)\pi}{m}+\frac{(2k-1)\pi}{m}\right).
\end{eqnarray*}
By setting $2k=m+1-2j$ for $S_{2}$ we have
$$
S_{2} = -\sum_{j =1}^{\frac{m-1}{2}} \E^{-(h\beta-h\gamma) \cos \frac{2\pi j}{m}}
\cos \left((h\gamma-h\beta) \sin \frac{2\pi j}{m}+\frac{2\pi j}{m}\right).
$$
Hence
$$
S_{2} = -\sum_{k = 1 }^{\frac{m-1}{2}} \E^{(h\gamma-h\beta) \cos \frac{2\pi k}{m}}
\cos \left((h\gamma-h\beta) \sin \frac{2\pi k}{m}+\frac{2\pi k}{m}\right).
$$
Now, using the formula $\cos(x+y)=\cos(x)\cos(y)-\sin(x)\sin(y)$, for $S_{1}$ and $S_{2}$
we get
\begin{eqnarray*}
S_{1} &=& \sum_{k = 1 }^{\frac{m-1}{2}}\Bigg[ \E^{-h\gamma \cos \frac{2\pi k}{m}}
\cos \left(h\gamma \sin \frac{2\pi k}{m}\right)   \E^{h\beta \cos \frac{2\pi k}{m}}
\cos \left(h\beta \sin \frac{2\pi k}{m}+\frac{2\pi k}{m}\right)\\
&&+\E^{-h\gamma \cos \frac{2\pi k}{m}}
\sin \left(h\gamma \sin \frac{2\pi k}{m}\right)  \E^{h\beta \cos \frac{2\pi k}{m}}
\sin \left(h\beta \sin \frac{2\pi k}{m}+\frac{2\pi k}{m}\right) \Bigg],\\
S_{2} &=& \sum_{k = 1 }^{\frac{m-1}{2}}\Bigg[ \E^{-h\beta \cos \frac{2\pi k}{m}}
\cos \left(h\beta \sin \frac{2\pi k}{m}\right)  \E^{h\gamma \cos \frac{2\pi k}{m}}
\cos \left(h\gamma \sin \frac{2\pi k}{m}+\frac{2\pi k}{m}\right)\\
&&+\E^{-h\beta \cos \frac{2\pi k}{m}}
\sin \left(h\beta \sin \frac{2\pi k}{m}\right)  \E^{h\gamma \cos \frac{2\pi k}{m}}
\sin \left(h\gamma \sin \frac{2\pi k}{m}+\frac{2\pi k}{m}\right) \Bigg].
\end{eqnarray*}
Substituting the sum (\ref{(66)}) of the last expressions for $S_{1}$ and $S_{2}$ into (\ref{(65)})
and using equalities (\ref{(31)})-(\ref{(34)}), after some simplifications, for $\beta<0$ we reduce the expression (\ref{(65)}) to the form
$$
v_m(h\beta)=-\frac{1}{2m}Q(h\beta)-P(h\beta,b_0, b_{1,k},b_{2,k}),
$$
where
\begin{eqnarray}
&&Q(h\beta)=\frac12 \E^{h\beta}(1-\E^{-1})\nonumber\\
&&\qquad+\sum_{k = 1 }^{\frac{m-1}{2}} \E^{h\beta \cos \frac{2\pi k}{m}} \Bigg[ \cos \left( h\beta \sin \frac {2\pi k}{m} \right) - \E^{-\cos \frac{2\pi k}{m}}
\cos \left(( h\beta-1) \sin \frac {2\pi k}{m} \right) \Bigg],\nonumber\\
&&P(h\beta, b_0, b_{1,k}, b_{2,k})=b_0 \E^{-h\beta}\nonumber\\
&&\qquad+\sum_{k = 1 }^{\frac{m-1}{2}} \E^{-h\beta \cos \frac{2\pi k}{m}}
\Bigg[ b_{1,k}\cos \left( h\beta \sin \frac {2\pi k}{m} \right) +b_{2,k}\sin \left( h\beta \sin \frac {2\pi k}{m} \right) \Bigg],\label{(67-1)}
\end{eqnarray}
where $b_0, b_{1,k}, b_{2,k},$ $k=1,2,...,\frac{m-1}{2}$ are unknowns.

Direct calculations show that $v_m(h\beta)$ for $\beta>N$ has the form
$$
v_m(h\beta)=\frac{1}{2m}Q(h\beta)+P(h\beta, b_0, b_{1,k}, b_{2,k}).
$$
Using last two expressions for $v_m(h\beta)$ and keeping in mind (\ref{(63)}) we have
\begin{equation}\label{(67-2)}
u_m(h\beta)=\left\{
\begin{array}{ll}
-\frac {1}{2m}Q(h\beta)+P(h\beta, d_0, d_{1,k}, d_{2,k}) - P(h\beta, b_0, b_{1,k}, b_{2,k}),& \beta<0,\\[2 mm]
\frac {1}{2m}Q(h\beta)+P(h\beta, d_0, d_{1,k}, d_{2,k}) + P(h\beta, b_0, b_{1,k}, b_{2,k}),  & \beta>N,\\
\end{array} \right.
\end{equation}
where $P(h\beta, d_0, d_{1,k}, d_{2,k})$ and $P(h\beta, b_0, b_{1,k}, b_{2,k})$ are defined by (\ref{(27)}) and
(\ref{(67-1)}), respectively. Here $d_0, d_{1,k}, d_{2,k}$ and $b_0, b_{1,k}, b_{2,k}$ are unknowns.

We denote
\begin{equation}\label{(67-3)}
\begin{array}{lll}
d_0^-=d_0-b_0,& d_{1,k}^-=d_{1,k}-b_{1,k},& d_{2,k}^-=d_{2,k}-b_{2,k},\ k=1,2,...,\frac{m-1}{2},\\[2mm]
d_0^+=d_0+b_0,& d_{1,k}^+=d_{1,k}+b_{1,k},& d_{2,k}^+=d_{2,k}+b_{2,k},\ k=1,2,...,\frac{m-1}{2}.
\end{array}
\end{equation}
If we find unknowns $d_0^-, d_{1,k}^-, d_{2,k}^-$ and $d_0^+, d_{1,k}^+, d_{2,k}^+$, $k=1,2,...,\frac{m-1}{2}$ then from (\ref{(67-3)})
we obtain $d_0, d_{1,k}, d_{2,k}$ and  $b_0, b_{1,k}, b_{2,k}$, $k=1,2,...,\frac{m-1}{2}$ as follows
\begin{equation}\label{(67-4)}
\begin{array}{lll}
d_0=\frac 12(d_0^++d_0^-),& d_{1,k}=\frac 12(d_{1,k}^++d_{1,k}^-),& d_{2,k}=\frac 12(d_{2,k}^++d_{2,k}^-),\ k=1,2,...,\frac{m-1}{2},\\[2mm]
b_0=\frac 12(d_0^+-d_0^-),& b_{1,k}=\frac 12(d_{1,k}^+-d_{1,k}^-),& b_{2,k}=\frac 12(d_{2,k}^+-d_{2,k}^-),\ k=1,2,...,\frac{m-1}{2}.
\end{array}
\end{equation}

Since for $\beta<0$ and $\beta>N$ the coefficients $C_{\beta}=0$ then from (\ref{(64)}) for $d_0^-, d_{1,k}^-, d_{2,k}^-$ and $d_0^+, d_{1,k}^+, d_{2,k}^+$,
$k=1,2,...,\frac{m-1}{2}$
we get the following system of $2m$ linear equations
$$
D_m(h\beta)*u_m(h\beta)=0\mbox{ for } \beta=-1,-2,...,-m \mbox{ and } \beta=N+1,N+2,...,N+m.
$$
Solving the last system, taking (\ref{(67-3)}) and (\ref{(67-4)}) into account, from equalities  (\ref{(31)}) and (\ref{(67-2)})
we get the following explicit form for $u_m(h\beta)$:
\begin{equation}\label{(67-5)}
u_m(h\beta)=\left\{
\begin{array}{ll}
-\frac {1}{2m}Q(h\beta)+P(h\beta, d_0^-, d_{1,k}^-, d_{2,k}^-),& \beta<0,\\[2 mm]
f_m(h\beta),& \beta=0,1,...,N,\\[2mm]
\frac {1}{2m}Q(h\beta)+P(h\beta, d_0^+, d_{1,k}^+, d_{2,k}^+),  & \beta>N,\\
\end{array} \right.
\end{equation}
Then we get coefficients $C_\beta$ of optimal quadrature formulas of the form (\ref{(1)}) as follows
$$
C_\beta=D_m(h\beta)*u_m(h\beta), \  \beta=0,1,...,N.
$$

Thus, Sard's problem on construction of optimal quadrature formulas of the form (\ref{(1)}) in the space $W_2^{(m,0)}$ for odd natural numbers $m$ is solved.

\begin{remark} It should be noted that for $m=1$, by realizing the above given algorithm, we get the optimal quadrature formula of the form (\ref{(1)}) in the space $W_2^{(1,0)}$ which was constructed in the works \cite{ShadHay04,ShadHay14}.
\end{remark}

\begin{remark} We note that using the discrete analogue $D_1(h\beta)$ of the operator $\frac{\d^2}{\d x^2}-1$ which corresponds to the case $m=1$ we
get interpolation spline $S_1(x)$ minimizing semi-norm in the space $W_2^{(1,0)}$ obtained in Theorem 3.2 of \cite{ShadHay13}.
It should be remembered that the spline $S_1(x)$ was used in \cite{AvezHay10} to
determine the total incident solar radiation at each time step and the temperature
dependence of the thermophysical properties of the air and water.
\end{remark}

\begin{remark} It is important to note that the discrete analogue $D_m(h\beta)$ of the differential operator $\frac{\d^{2m}}{\d x^{2m}}-1$, constructed
in Sect. 5, can be applied for construction of interpolation splines minimizing the semi-norm (\ref{norm}) and optimal quadrature formulas for
numerical integration of Fourier coefficients in the space $W_2^{(m,0)}$.
\end{remark}

In the next section we give the results of the algorithm for construction of optimal quadrature formulas, described in this section, for $m=1$ and $m=3$.

\section{Coefficients of the optimal quadrature formulas for $m=1$ and $m=3$}

In this section we give the results of the realization of the algorithm for construction of optimal quadrature
formulas (\ref{(1)}) in the space $W_2^{(m,0)}$ for the cases $m=1$ and $m=2$.

We recall that in the case $m=1$ we get the result of Theorem 4 of \cite{ShadHay04} and Theorem 4.4 of \cite{ShadHay14}.
For $m=1$ the system (\ref{(31)})-(\ref{(34)}) has the form
\begin{eqnarray}
&&C_{\beta}*G_1(h\beta)+d_0\E^{-h\beta}=f_1(h\beta),\ \beta=0,1,...,N,\label{(73)}\\
&&\sum_{\beta  = 0}^N C_\beta  \E^{-h \beta}  = 1 -\frac {1}{\E},\label{(74)}
\end{eqnarray}
where
\begin{eqnarray*}
G_1(h\beta)&=&\frac{\mathrm{sgn} x}{2}\sinh(x),\\
f_1(h\beta)&=&\frac{1}{2}\left[\cosh(h\beta)+\cosh(1-h\beta)\right]-1,
\end{eqnarray*}

The system (\ref{(73)})-(\ref{(74)}) was solved in the work \cite{ShadHay04} and the following theorem was proved.

\begin{theorem} (Theorem 4 of \cite{ShadHay04}). The coefficients of optimal quadrature formulas
of the form (\ref{(1)}) with equally spaced nodes in the space $W_2^{(1,0)}$ are expressed by
formulas
\begin{eqnarray*}
C_0&=&\frac{\E^h-1}{\E^h+1},\\
C_{\beta}&=&\frac{2(\E^h-1)}{\E^h+1},\ \beta=1,2,...,N-1,\\
C_N&=&\frac{\E^h-1}{\E^h+1},
\end{eqnarray*}
where $h=1/N,$ $N=1,2,...$.
\end{theorem}

\begin{remark} It was shown in Sect. 5 of \cite{ShadHay14} that the error of the
optimal quadrature formula in the space $W_1^{(1,0)}$ is less than the error of the optimal
quadrature formula in the space $L_2^{(1)}$.
\end{remark}

For $m=3$ the system (\ref{(31)})-(\ref{(34)}) has the form
\begin{eqnarray}
&&C_\beta*G_3(h\beta) +  d_0 \E^{-h \beta} + d_{1,1} \E^{\frac {h \beta}{2}} \cos \left(\frac {\sqrt{3}}{2} h \beta\right)\nonumber \\
&&\qquad\qquad +d_{1,2} \E^{\frac {h \beta}{2}} \sin \left(\frac {\sqrt{3}}{2} h\beta \right)= f_3(h\beta), \ \ \beta  = 0,1,\ldots,N,\label{(68)}\\
&&\sum_{\beta  = 0}^N C_\beta  \E^{-h \beta}  = 1 -\frac {1}{\E},\label{(69)}\\
&&\sum_{\beta  = 0}^N C_\beta \E^{\frac {h\beta}2} \cos \left( \frac{\sqrt{3}}{2} h \beta\right) =
\cos\frac{2\pi}{3}-\E^{\cos\frac{\pi}{3}}\cos\left(\sin\frac{\pi}{3}+\frac{2\pi}{3}\right), \label{(70)}\\
&&\sum_{\beta  = 0}^N C_\beta \E^{\frac {h \beta}2} \sin \left( \frac{\sqrt{3}}{2} h \beta\right) =
\sin\frac{2\pi}{3}-\E^{\cos\frac{\pi}{3}}\sin\left(\sin\frac{\pi}{3}+\frac{2\pi}{3}\right),\label{(71)}
\end{eqnarray}
where
\begin{eqnarray*}
G_3(x)&=&\frac{\mathrm{sgn}x}{6} \left(\sinh(x)+\sum_{n=1}^{2}\E^{x  \cos\left(\frac{\pi n}{3}\right)}
\cos\left(x\sin\left(\frac{\pi n}{3}\right)+\frac{\pi n}{3}\right)\right),\\
f_3(h\beta)&=& \frac{1}{12}\left[\E^{-h\beta}(\E+1)+\E^{h\beta}(\E^{-1}+1)+2\E^{-\frac{h \beta}{2}} \cos \left(\frac {\sqrt{3}}{2}h \beta \right) \right.\left(\E^{\frac 12}\cos\left( \frac {\sqrt{3}}{2}\right)+1\right) \nonumber\\
&& +2\E^{\frac 12}
\sin \left( \frac {\sqrt{3}}{2}\right)\E^{-\frac{h \beta}{2}} \sin \left( \frac{\sqrt{3}}{2} h \beta\right)
+2\E^{-\frac 12} \sin \left(\frac {\sqrt{3}}{2}\right) \E^{\frac{h \beta}{2}} \sin \left( \frac{\sqrt{3}}{2} h \beta\right)\nonumber \\
&&\left.+2\E^{\frac{h \beta}{2}} \cos \left(\frac {\sqrt{3}}{2}h \beta \right)
\cdot \left(\E^{-\frac 12}\cos\left( \frac {\sqrt{3}}{2}\right)+1\right)-1\right].
\end{eqnarray*}

The system (\ref{(68)})-(\ref{(71)}) was solved in \cite{Boltaev14} and the following result was obtained.

\begin{theorem} The coefficients of the optimal quadrature formula (\ref{(1)}) in the space
$W_2^{(3,0)}$ have the form
\begin{eqnarray*}
C_0& =& 1 - \left( {\frac{T}{{{\E^h} - 1}} + \frac{{{m_1}{\tau _1}}}{{{\E^h} - {\tau _1}}} +
\frac{{{m_2}{\tau _2}}}{{{\E^h} - {\tau _2}}} + \frac{{{n_1}\tau _1^N}}{{{\tau _1}{\E^h} - 1}} +
\frac{{{n_2}\tau _2^N}}{{{\tau _2}{\E^h} - 1}}} \right),\\
C_{\beta } &=& T + {m_1}\tau _1^\beta  + {m_2}\tau _2^\beta  +
{n_1}\tau _1^{N - \beta } + {n_2}\tau _2^{N - \beta },\,\,\,\,\beta  = 1,...,N - 1,\\
C_N& =&  - 1 + {\E^h}\left( {\frac{T}{{{\E^h} - 1}} + \frac{{{m_1}\tau _1^N}}{{{\E^h} -
{\tau _1}}} + \frac{{{m_2}\tau _2^N}}{{{\E^h} - {\tau _2}}} +
\frac{{{n_1}{\tau _1}}}{{{\tau _1}{\E^h} - 1}} +
\frac{{{n_2}{\tau _2}}}{{{\tau _2}{\E^h} - 1}}} \right),
\end{eqnarray*}
where $m_1,\ m_2,\ n_1$ and $n_2$ satisfy the following system
\begin{eqnarray*}
&&A_{11}m_1+A_{12}m_2+\tau_1^NB_{11}n_1+\tau_2^NB_{12}n_2=T_1,\\
&&A_{21}m_1+A_{22}m_2+\tau_1^NB_{21}n_1+\tau_2^NB_{22}n_2=T_2,\\
&&\tau_1^NA_{11}m_1+\tau_2^NA_{12}m_2+B_{11}n_1+B_{12}n_2=T_1,\\
&&\tau_1^NA_{21}m_1+\tau_2^NA_{22}m_2+B_{21}n_1+B_{22}n_2=T_2,
\end{eqnarray*}
here $$
\begin{array}{ll}
\displaystyle A_{1k}=\frac{\tau_k\E^{h/2}\sin(\sqrt{3}h/2)}{1-2\tau_k\E^{h/2}\cos(\sqrt{3}h/2)+\tau_k^2\E^h},&
\displaystyle A_{2k}=\frac{\E^h}{\E^h-\tau_k}+\frac{\tau_k\E^{h/2}\cos(\sqrt{3}h/2)-1}{1-2\tau_k\E^{h/2}\cos(\sqrt{3}h/2)+\tau_k^2\E^h},\\[5mm]
\displaystyle B_{1k}=\frac{\tau_k\E^{h/2}\sin(\sqrt{3}h/2)}{\tau_k^2-2\tau_k\E^{h/2}\cos(\sqrt{3}h/2)+\E^h},&
\displaystyle B_{2k}=\frac{\E^h\tau_k}{\E^h\tau_k-1}+\frac{\tau_k\E^{h/2}\cos(\sqrt{3}h/2)-\tau_k^2}{\tau_k^2-2\tau_k\E^{h/2}\cos(\sqrt{3}h/2)+\E^h},\\[5mm]
\end{array}
$$
for $k=1,2$ and
$$
T_1=\frac{\sqrt{3}}{2}-\frac{T\E^{h/2}\sin(\sqrt{3}h/2)}{1-2\E^{h/2}\cos(\sqrt{3}h/2)+\E^h},\
T_2=\frac{3}{2}-\frac{T\E^h}{\E^h-1}-\frac{T\E^{h/2}\cos(\sqrt{3}h/2)-T}{1-2\E^{h/2}\cos(\sqrt{3}h/2)+\E^h},
$$
\begin{eqnarray*}
\tau _1& =& \frac{1}{4}\left[ {{K_1} + \sqrt {K_1^2 - 4{K_2} + 8}  + \sqrt {{{\left( {{K_1} +
    \sqrt {K_1^2 - 4{K_2} + 8} } \right)}^2} - 16} } \right],\\
\tau _2& =& \frac{1}{4}\left[ {{K_1} + \sqrt {K_1^2 - 4{K_2} + 8}  - \sqrt {{{\left( {{K_1} +
    \sqrt {K_1^2 - 4{K_2} + 8} } \right)}^2} - 16} } \right],\\
T &=&\frac{{{24\left( {\cosh(h)-1}\right){{\left({\cos
\left({\frac{{\sqrt 3 }}{2}h} \right) - \cosh\left( {\frac{h}{2}}
\right)} \right)}^2}}}}{{K\left( {{K_2} + 2 - 2{K_1}} \right)}},\\
K &=& {\sinh(h)+\sinh\left({\frac{h}{2}}\right) \cos \left( {\frac{{\sqrt 3 }}{2}h} \right) -
\sqrt 3 \cosh\left({\frac{h}{2}}\right)\sin \left({\frac{{\sqrt 3 }}{2}h} \right)},\\
K_1 &=& {2\cosh(h)}+ {\frac{{4\cos \left( {\frac{{\sqrt 3 }}{2}h} \right)\cosh\left({\frac{h}{2}}\right)\sinh(h) +
\sinh(h) - \sqrt 3 \sin \left({\sqrt 3 h} \right) - 2\sinh(h)\cosh(h)}}{{\sinh(h) +
\sinh\left({\frac{h}{2}}\right) \cos \left({\frac{{\sqrt 3 }}{2}h} \right) - \sqrt 3 \cosh\left( {\frac{h}{2}} \right)
 \sin \left( {\frac{{\sqrt 3 }}{2}h} \right)}}},\\
K_2& =& {\frac{{2\cos \left( {\sqrt 3 h} \right)\sinh(h) + 4\sinh(h)\cosh(h) - 2\sqrt 3 \sin \left( {\sqrt 3 h} \right)\cosh(h)}}
{{\sinh(h)+\sinh\left({\frac{h}{2}}\right) \cos \left( {\frac{{\sqrt 3 }}{2}h} \right) -
\sqrt 3 \cosh\left( {\frac{h}{2}} \right) \sin \left( {\frac{{\sqrt 3 }}{2}h} \right)}}}.
\end{eqnarray*}
\end{theorem}

\section*{Acknowledgments}

This work has been done while A.R.Hayotov was visiting Department of Mathematical
Sciences at KAIST, Daejeon, Republic of Korea. A.R.Hayotov is very grateful to professor Chang-Ock Lee and his research group for hospitality.
A.R. Hayotov's work was supported by the 'Korea Foundation for Advanced Studies'/'Chey Institute for Advanced Studies' International Scholar Exchange Fellowship for academic year of 2018-2019

\end{document}